\documentclass[12pt]{amsart}
\usepackage{amscd,amsmath,amssymb,amsfonts}
\usepackage[cmtip, all]{xy}
\theoremstyle{plain}
\newtheorem{thm}{Theorem}

\newtheorem{cor}[thm]{Corollary}

\theoremstyle{definition}

\newtheorem{rmk}[thm]{Remark}

\numberwithin{thm}{section}
\numberwithin{equation}{section}

\newcommand{\ml}[2]{\begin{multline}\label{#1}#2 \end{multline}}
\newcommand{\ga}[2]{\begin{gather}\label{#1}#2 \end{gather}}

\newcommand{\sI}{{\mathcal I}}

\newcommand{\sO}{{\mathcal O}}

\newcommand{\A}{{\mathbb A}}

\renewcommand{\P}{{\mathbb P}}
\newcommand{\Q}{{\mathbb Q}}

\begin{document}

\title[Appendix]{Appendix to ``Congruences for rational points on varieties
over finite fields'' by N. Fakhruddin and C. S. Rajan
 }
\author{H\'el\`ene Esnault}
\address{
Universit\"at Essen, FB6, Mathematik, 45117 Essen, Germany}
\email{esnault@uni-essen.de}

\date{March 12, 2004}
\maketitle
We give a positive answer to questions 4.2 and 4.4 of \cite{FR}.
\section{Hodge type}
In this section we consider varieties $X$ defined over a field of characteristic 0, their de Rham cohomology $H^m(X)$ together with their Hodge filtration $F$. 
\begin{thm} \label{thm:HT} Let $f: X\to Y$ be a proper surjective
morphism between smooth varieties defined over a field of characteristic 0
with $Y$ connected. If there is a closed point $y\in Y$ such that $f$ is smooth in $y$ and $gr_0^FH^m(f^{-1}(y))=0$ for all $m\ge 1$, then 
$gr_0^FH^m(f^{-1}(z)_{{\rm red}})=0$ for all closed points $z\in Y$ and all $m\ge 1$. \end{thm}
\begin{proof}
Since $f^{-1}(y)$ is smooth, $gr_0^FH^m(f^{-1}(y))=H^m(f^{-1}(y), \sO)=0$
for all $m\ge 1$.  
Since $R^mf_*\sO_X$ is locally free on the smooth locus of $f$, then $R^mf_*\sO_X$ has to be a torsion sheaf on $Y$. Since $X$ and $Y$ are smooth and $f$ is surjective, by Koll\'ar's torsion-freeness 
theorem \cite{Ko}, one has
\ga{1}{ R^mf_*\sO_X=0.}
 Let $\frak{m}_z$ be the maximal ideal of a closed point $z$. Since $Y$ is smooth, locally on $Y$, there are smooth divisors $D_1,\ldots, D_a$, $a={\rm dim}(Y)$ with $$\frak{m}_z=
(\sO_Y(-D_1),\ldots, \sO_Y(-D_a)).$$ Thus one has a resolution
\ml{2}{0\to \sO_X(-\sum_1^a f^*(D_i))\to \ldots \to \oplus_{i_1<\ldots < i_r}
\sO_X(-f^*(D_{i_1}+\ldots + D_{i_r}))\\
\to \ldots \to \sO_X \to \sO_{f^{-1}(z)}\to 0.}
Thus the vanishing \eqref{1} and projection formula imply 
\ga{3}{H^m(X, \sO_{f^{-1}(z)})=0 \ \forall \ m \ge 1.}
Let now $f': X'\to Y'$ be a compactification of $f$, that is $X', Y'$ are smooth
 proper, $X\subset X'$ and $Y\subset Y'$ are open dense, and $f'=f|_Y$.
Let
 $\sI=f^{'*}\frak{m}_z$ be the ideal sheaf of the fiber $f^{'-1}(z)$. 
Then by \cite{E}, Proposition 1.2, the  map $\sigma^*$
\ga{4}{\sigma^*: H^n(X', \sI)\to 
gr_0^FH^n(\tilde{X}, (\sigma^{-1}f^{'-1}(z))_{{\rm red}})=
gr_0^FH^n(X', (f^{'-1}(z))_{{\rm red}})}
is surjective for all $n$, 
 where $\sigma: \tilde{X}\to X'$ is a birational map, isomorphic over $X'\setminus f^{'-1}(z)$, 
 with $\tilde{X}$ smooth,  and $\sigma^{-1}f^{'-1}(z)$ 
a normal crossing divisor. 
From the commutative diagram of exact sequences 
\ga{5}{\begin{CD}
H^m(X', \sO_{X'}) @>>> H^m(X', \sO_{f^{'-1}(z)}) @>>> H^{m+1}(X', \sI)\\
@V= VV  @V\sigma^* VV @VV\sigma^* V \\
gr_0^F H^m(X') @>>> gr_0^F H^m(f^{-1}(z)_{{\rm red}}) @>>> 
gr_0^FH^n(X', (f^{'-1}(z))_{{\rm red}})
\end{CD}
}
we conclude that the vertical map in the middle is surjective, thus by \eqref{3} that $gr_0^F H^m(f^{-1}(z)_{{\rm red}})=0$. 
\end{proof}
\begin{rmk} Because of the resolution \eqref{2} it is not possible to separate the $m$ involved in the proof of the Theorem.  
\end{rmk} 
\begin{cor} \label{cor}
Let $f: X\to Y$ be a proper morphism between smooth varieties defined over a field $k$ of characteristic 0 with $Y$ connected. If there is a closed point $y\in Y$
such that f is smooth in $y$ and $CH_0(f^{-1}(y)\times_k \overline{k(f^{-1}(y))})=\Q$, then 
$gr_0^FH^m(f^{-1}(z)_{{\rm red}})=0$ for all closed points $z\in Y$ and all $m\ge 1$. \end{cor}
\begin{proof}
By S. Bloch's theorem, \cite{B}, Appendix to Lecture 1, the assumption implies
$H^m(f^{-1}(y), \sO)=gr_0^FH^m(f^{-1}(y))=0$ for all $m\ge 1$. One applies Theorem \ref{thm:HT}.
\end{proof}
\begin{rmk} 
According to Bloch's conjecture, the assumptions in Theorem 
\ref{thm:HT} and in its corollary should be equivalent. As we are far from knowing Bloch's conjecture, we have formulated Theorem \ref{thm:HT} purely in the coherent category.
\end{rmk} 

\section{Eigenvalues of Frobenius}
In this section, we consider varieties $X$ defined over a finite field $k$, their \'etale cohomology $H^m(\overline{X}, \Q_\ell)$,  with $\overline{X}=X\times_k \overline{k}$, and the eigenvalues of the geometric Frobenius $F$  acting on it.  
\begin{thm} \label{thm:EV} 
Let $f: X\to Y$ be a proper surjective morphism  of smooth irreducible varieties defined over a finite field $k$ with $q$ elements. Assume $CH_0((X\times_Y \overline{k(X))})=\Q$. Then for all closed points $z\in Y(k)$, the eigenvalues of the geometric
Frobenius  $F$ acting on $H^m(\overline{f^{-1}(y)_{{\rm red}}}, \Q_\ell)$ are divisible by $q$ as algebraic integers for all $m\ge 1$. 
\end{thm}
\begin{proof}
We may assume that $Y$ is affine of dimension $a$, 
containing $y$. We set $Y'=Y\setminus \{y\}, 
X'=X\setminus \{f^{-1}(y)_{{\rm red}}\}$.  
Then 
\ga{2.1}{H^n_c(\overline{Y'}, \Q_\ell)\to H^n_c(\overline{Y}, \Q_\ell)} 
is surjective for $n\ge 1$ and an isomorphism for $n\ge 2$ .
Thus via the exact sequence
\ml{2.2}{\ldots \to H^n_c(\overline{X'}, \Q_\ell) 
\to H^n_c(\overline{X}, \Q_\ell) \to \\ 
H^n(\overline{f^{-1}(y)_{{\rm red}}}, \Q_\ell) \to H^{n+1}_c(\overline{X'}, \Q_\ell)\to \ldots}
one concludes that the image of $H^n_c(\overline{X}, \Q_\ell)$ in 
$H^n 
(\overline{f^{-1}(y)_{{\rm red}}}, \Q_\ell)$ is a quotient of 
$H^n_c(\overline{X}, \Q_\ell)/H^n_c(\overline{Y}, \Q_\ell)$
thus by 3.1, Proof of Theorem 1.1 of \cite{FR}, the Frobenius eigenvalues on it are divisible by $q$. Thus we just have to consider
the image of 
$H^n 
(\overline{f^{-1}(y)_{{\rm red}}}, \Q_\ell)$ in
$H^{n+1}_c(\overline{X'}, \Q_\ell)$ for $n\ge 1$. 
Again by the isomorphism  \eqref{2.1}  for $n\ge 2$,  
 by Deligne's integrality theorem 
\cite{DeInt}, Corollaire 5.5.3, together with 
Artin's vanishing theorem \cite{Ar}, Th\'eor\`eme 3.1, and the
 Proof of Theorem 1.1 of \cite{FR},
the only $n\ge 1$ for which possibly the Frobenius eigenvalues on 
$H^{n+1}_c(\overline{X'}, \Q_\ell)$, and consequently on
$H^n 
(\overline{f^{-1}(y)_{{\rm red}}}, \Q_\ell)$ 
 are not divisible by $q$ is $n=a-1$.  
But $f^{-1}(y)_{{\rm red}}$ is geometrically irreductible.
The Lefschetz trace formula implies then that the eigenvalues of Frobenius on
 $H^{a-1} 
(\overline{f^{-1}(y)_{{\rm red}}}, \Q_\ell)$ are divisible by $q$ if and only
if
 $f^{-1}(y)_{{\rm red}}$ has one $k$-point modulo $q$. This 
is  Theorem 1.1 of \cite{FR}.
 
\end{proof}
\begin{rmk}
 One can avoid the Euler characteristic argument via the Lefschetz trace formula in the proof of Theorem \ref{thm:EV} by replacing $f:X \to Y$ by $
f\times {\rm id}: X\times \A^1\to Y\times \A^1$ for example, as observed by N. Fakhruddin. However, one has to say that the version of Bloch's decomposition of the diagonal used in 3.1, Proof of Theorem 1.1 of \cite{FR} is equivalent
to the assumption of Theorem 1.1 so as to see that one may replace $\overline{k(X)}$ by $\overline{k(X\times \A^1)}$ in the assumption. 
This gives of course as well a motivic proof of Corollary \ref{cor}, however not of Theorem \ref{thm:HT} in absence of a solution to Bloch's conjecture.  
\end{rmk}

\bibliographystyle{plain}

\renewcommand\refname{References}

\end{document}